\documentclass[12pt]{amsart}
\usepackage{amsmath}
\theoremstyle{plain}
\newtheorem{thm}{Theorem}[section]
\newtheorem{prop}[thm]{Proposition}
\newtheorem{lemma}[thm]{Lemma}
\newtheorem{cor}[thm]{Corollary}

\theoremstyle{definition}
\newtheorem{example}{Example}

\newcommand{\C}{{\mathbb C}}
\newcommand{\Z}{{\mathbb Z}}

\begin{document}
\thispagestyle{empty}

\title{A Note on the Symmetric Powers of the Standard Representation
  of $S_n$} 
\maketitle

\begin{center}
{
David Savitt\footnote{Supported by an NSERC PGS-B fellowship} \\
Department of Mathematics, Harvard University \\
Cambridge, MA 02138, USA \\
{\small\texttt{dsavitt@math.harvard.edu}} \\

\hskip 0.5cm

Richard P. Stanley\footnote{Partially supported by
  NSF grant 
  DMS-9500714} \\
Department of Mathematics, Massachusetts Institute of Technology \\
Cambridge, MA 02139, USA \\
{\small\texttt{rstan@math.mit.edu}}
} \\

\hskip 0.5cm

Submitted: January 7, 2000; Accepted: February 12, 2000

\end{center}

\section*{Abstract}

\pagestyle{myheadings} \markboth{\hfill{\sc the electronic journal of
combinatorics 7 (2000), \#Rxx}} {{\sc the
              electronic journal of combinatorics 7 (2000), \#Rxx}\hfill}

In this paper, we prove that the dimension of the space
spanned by the characters of the symmetric powers of the
standard $n$-dimensional representation of $S_n$ is
asymptotic to $n^2 / 2$.  This is proved by using generating
functions to obtain formulas for upper and lower bounds,
both asymptotic to $n^2/2$, for this dimension.  In
particular, for $n \ge 7$, these characters do not span the
full space of class functions on $S_n$.
 
\section*{Notation}

Let $P(n)$ denote the number of (unordered) partitions of $n$ into
positive integers, and let $\phi$ denote the Euler totient function.  Let 
$V$ be the standard $n$-dimensional representation
of $S_n$, so that $V = \C e_1 \oplus \cdots \oplus \C e_n$ with
$\sigma(e_i)=e_{\sigma i}$ for $\sigma \in S_n$.  Let $S^N V$ denote the
$N^{\rm th}$ symmetric power of $V$, and let $\chi_N : S_n \rightarrow \Z$
denote its character.  Finally, let $D(n)$ denote the dimension of the
space of class functions on $S_n$ spanned by all the $\chi_N$, $N \ge 0$. 

\section{Preliminaries}

Our aim in this paper is to investigate the numbers $D(n)$.  It is a
fundamental problem of invariant theory to decompose the character of the
symmetric powers of an irreducible representation of a finite group (or
more generally a reductive group). A special case with a nice theory is
the reflection representation of a finite Coxeter group. This is
essentially what we are looking at. (The defining representation of $S_n$
consists of the direct sum of the reflection representation and the
trivial representation. This trivial summand has no significant effect on
the theory.) In this context it seems natural to ask: what is the
dimension of the space spanned by the symmetric powers? Moreover,
decomposing
the symmetric powers of the character of an irreducible representation of
$S_n$ is an example of the operation of \emph{inner plethysm} \cite[Exer.\
7.74]{EC2}, so we are also obtaining some new information related to this
operation.

We begin with:

\begin{lemma} 

 Let $\lambda = (\lambda_1,\dots,\lambda_k)$ be a partition of $n$
 (which we  
 denote by $\lambda\vdash n$), and suppose $\sigma \in S_n$ is a
 $\lambda$-cycle.  Then $\chi_N (\sigma)$ is equal to the number of
 solutions 
 $(x_1,\dots,x_k)$ in nonnegative integers to the equation $\lambda_1 x_1 +
 \cdots + \lambda_k x_k = N$.

\end{lemma}

\begin{proof} 

 Suppose without loss of generality that $\sigma = (1 \ 2 \
 \cdots \ \lambda_1) (\lambda_1 + 1 \ \cdots \ \lambda_1 + \lambda_2) 
 \cdots ( \lambda_1 + \cdots + \lambda_{k-1} + 1 \ \cdots \
 n)$. Consider a  
 basis vector $e_1^{\otimes c_1} \otimes \cdots \otimes e_n^{\otimes
 c_n}$ of 
 $S^N V$, so that $c_1 + \cdots + c_n = N$ with each $c_i \ge 0$.  This
vector is fixed by
 $\sigma$ 
 if and only if $c_1 = \cdots = c_{\lambda_1}$, $c_{\lambda_1 + 1} =
 \cdots = 
 c_{\lambda_1 + \lambda_2}$ and so on.  Since $\chi_N (\sigma)$ equals the
 number of basis vectors fixed by $\sigma$, the lemma follows. \end{proof}

It seems difficult to work directly with the $\chi_N$'s; fortunately,
it is not too hard to restate the problem in more concrete terms.
 Given a partition $\lambda = (\lambda_1,\dots,\lambda_k)$ of $n$,  define

 \begin{equation} f_\lambda(q) =
   \frac{1}{\left(1-q^{\lambda_1}\right)\cdots 
   \left(1-q^{\lambda_k}\right)}. \label{eq:f} 
 \end{equation}

Next, define $ F_n \subset \C[[q]] $ to be the complex vector space
spanned by all of these $f_{\lambda}(q)$'s.  We have:

\begin{prop}  \label{prop:dfn}
 $\dim F_n = D(n)$.
\end{prop}

\begin{proof}

 Consider the table of the characters $\chi_N$; we are interested in the
 dimension of the row-span of this table.  Since the dimension of the row-span
 of a matrix is equal to the dimension of its column-span, we can equally
 well study the dimension of the space spanned by the columns of the table.
 By the preceeding lemma, the $N^{\rm th}$ entry of the column
corresponding to
 the $\lambda$-cycles is equal to the number of nonnegative integer
 solutions to the equation $\lambda_1 x_1 + \cdots + \lambda_k x_k = N$. 
 Consequently,  one easily verifies that $f_{\lambda}(q)$ is the generating
 function for the entries of the column corresponding to the $\lambda$-cycles.  
 The dimension of the column-span of our table is therefore equal to $\dim 
 F_n$, and the proposition is proved.

\end{proof}

\section{Upper Bounds on $D(n)$}

Our basic strategy for computing upper bounds for $\dim F_n$ is to put all
the generating functions $f_{\lambda}(q)$ over a common denominator;
then the dimension of their span is bounded above by $1$ plus the degree
of their numerators.  For example, one can see without much difficulty
that $(1-q)(1-q^2)\cdots(1-q^n)$ is the least common multiple of the 
denominators of the $f_{\lambda}(q)$'s.  Putting all of the 
$f_{\lambda}(q)$'s
over this common denominator, their numerators then have degree 
$n(n+1)/2 - n$, which proves 
   \begin{equation} D(n) \le \frac{n(n-1)}{2} + 1.
       \label{eq:dn} \end{equation}

By modifying this 
strategy carefully, it is possible to find a somewhat better bound.
Observe that the denominator of each of our $f_{\lambda}$'s is (up to sign 
change) 
a product of cyclotomic polynomials.  In fact, the power of the
$j^{\rm th}$
cyclotomic polynomial $\Phi_j(q)$ dividing the denominator of
$f_{\lambda}(q)$ is precisely equal to the number of $\lambda_i$'s
which are 
divisible by $j$.   It follows that $\Phi_j(q)$ divides the denominator
of $f_{\lambda}(q)$ at most $\left\lfloor \frac{n}{j} \right\rfloor$ times,
and the partitions $\lambda$ for which this upper bound is achieved
are precisely the $P\left(n - j \left\lfloor \frac{n}{j}
  \right\rfloor\right)$ partitions of $n$ 
which contain $\left\lfloor \frac{n}{j} \right\rfloor$ copies of $j$.
Let $S_j$ be 
the collection of $f_{\lambda}$'s corresponding to these $P\left(n - j
\left\lfloor  
\frac{n}{j} \right\rfloor\right)$ partitions.  One sees immediately
that the dimension of the space
spanned by the functions in $S_j$ is just $D\left(n - j \left\lfloor
    \frac{n}{j}  
\right\rfloor\right)$:  in fact, the functions in this space are exactly 
$1/(1-q^j)^{\left\lfloor \frac{n}{j} \right\rfloor}$ times the
functions in  
$F_{n - j \left\lfloor \frac{n}{j} \right\rfloor}$. 

Now the power of $\Phi_j(q)$ in the least common multiple of the 
denominators of all of the $f_{\lambda}(q)$'s {\em excluding those in $S_j$} 
is only $\left\lfloor \frac{n}{j} \right\rfloor - 1$, so the degree of this common
denominator is only $n(n+1)/2 - \phi(j)$.  Therefore, as in the
first paragraph of this section, the dimension of the space spanned by all 
of the $f_{\lambda}$'s except those in $S_j$ is at most $n(n-1)/2 + 1 -
\phi(j)$; since the dimension spanned by the functions in $S_j$ is 
$D\left(n - j \left\lfloor \frac{n}{j} \right\rfloor\right)$, we have
proved the upper bound
$$ D(n) \le \frac{n(n-1)}{2} + 1 - \phi(j) + D\left( n - j \left\lfloor
  \frac{n}{j}  
\right\rfloor\right) .$$  If it happens that $D\left(n - j
  \left\lfloor \frac{n}{j} \right\rfloor\right) < 
\phi(j)$, then this upper bound is an improvement on our original upper
bound.  If we repeat this process, this time simultaneously excluding the 
sets $S_j$ for {\em all} of the $j$'s which gave us an improved upper bound 
in the above argument, we find that we have proved:

\begin{prop} 
$$ D(n) \le \frac{n(n-1)}{2} + 1 - \sum_{j=1}^{n}
\max{\left(0 , \phi(j) - D\left(n - j \left\lfloor \frac{n}{j}
\right\rfloor\right)\right)}.$$
\end{prop}

Finally, we obtain an upper bound for $D(n)$ which does not depend on other
values of $D(\cdot)$:

\begin{cor} 
Recursively define $U(0) = 1$ and
$$ U(n) = \frac{n(n-1)}{2} + 1 - \sum_{j=1}^{n}
\max{\left(0 , \phi(j) - U\left(n - j \left\lfloor \frac{n}{j} \right\rfloor\right)\right)}.$$  
Then $D(n) \le U(n)$.
\end{cor}

\begin{proof}
We proceed by induction on $n$.  Equality certainly holds for $n = 0$.  For 
larger $n$,
the inductive hypothesis shows that $D\left(n - j \left\lfloor \frac{n}{j} \right\rfloor\right)
\le U\left(n - j \left\lfloor \frac{n}{j} \right\rfloor\right)$ when $j > 0$, and so
\begin{eqnarray*}
D(n) & \le & \frac{n(n-1)}{2} + 1 - \sum_{j=1}^{n}
\max{\left(0 , \phi(j) - D\left(n - j \left\lfloor \frac{n}{j} \right\rfloor\right)\right)} \\
      & \le & \frac{n(n-1)}{2} + 1 - \sum_{j=1}^{n}
\max{\left(0 , \phi(j) - U\left(n - j \left\lfloor \frac{n}{j} \right\rfloor\right)\right)} \\
      & =  & U(n).
\end{eqnarray*}
\end{proof}

Below is a table of values of $D(n)$ and $U(n)$ for $n \le 23$, calculated 
in Maple, with $P(n)$ and 
our first estimate $\frac{n(n-1)}{2} + 1$ provided for contrast. Note
that in the range $1\leq n\leq 23$, we have $D(n)=U(n)$ except for
$n=19,20$, when $U(n)-D(n)=1$. Is it true, for instance, that
$$-D(n) + \frac{n(n-1)}{2} + 1 - \sum_{j=1}^{n}
\max{\left(0 , \phi(j) - D\left(n - j \left\lfloor \frac{n}{j} 
\right\rfloor\right)\right)}$$ 
is bounded as
$n\rightarrow\infty$?

\begin{table}[h]
\begin{tabular}{|c|r|r|r|r|r|r|r|r|r|r|r|r|r|r|}
\hline
 $n$   & 1 & 2 & 3 & 4 & 5 & 6 & 7 & 8 & 9 & 10 & 11 & 12 & 13 & 14 \\ 
\hline
$D(n)$ & 1 & 2 & 3 & 5 & 7 & 11 & 13 & 19 & 23 & 29 & 35 & 45 & 51 & 62  \\
$U(n)$ & 1 & 2 & 3 & 5 & 7 & 11 & 13 & 19 & 23 & 29 & 35 & 45 & 51 & 62  \\ 
\hline
$n(n-1)/2 + 1$ & 1 & 2 & 4 & 7 & 11 & 16 & 22 & 29 & 37 & 46 & 56 & 67 & 
79 & 92  \\
$P(n)$ & 1 & 2 & 3 & 5 & 7 & 11 & 15 & 22 & 30 & 42 & 56 & 77 & 101 & 135 
 \\ \hline \end{tabular}

\vskip 0.1in

\begin{tabular}{|c|r|r|r|r|r|r|r|r|r|}
\hline
 $n$ &  15 & 16 & 17 & 18 & 19 & 20 & 21 & 22 & 23  \\ \hline
$D(n)$  & 69 & 79 & 90 & 106 & 118 & 134 & 146 & 161 & 176 \\
$U(n)$  & 69 & 79 & 90 & 106 & 119 & 135 & 146 & 161 & 176 \\ 
\hline
$n(n-1)/2 + 1$  & 106 & 121 & 137 & 154 & 172 & 191 & 211 & 232 & 254 \\
$P(n)$ & 176 & 231 & 297 & 385 & 490 & 627 & 792 & 1002 & 1255  \\ \hline
\end{tabular}

\vskip 0.1in

\caption{Values of $D(n)$, $U(n)$, $n(n-1)/2 + 1$, $P(n)$ for small $n$}
\end{table}

\begin{example}
The first dimension where $D(n) < P(n)$ is $n=7$, and it is easy then to
show that $D(n) < P(n)$ for all $n \ge 7$.  The difference $P(7)-D(7)=2$
arises from the following two relations:
$$ \frac{4}{(1-x^2)^2 (1-x)^3} = \frac{3}{(1-x^3)(1-x)^4} + 
\frac{1}{(1-x^3)(1-x^2)^2} $$ and
$$ \frac{3}{(1-x^3)(1-x^2)(1-x)^2} = \frac{2}{(1-x^4)(1-x)^3} + 
\frac{1}{(1-x^4)(1-x^3)}.$$

The first relation, for example, says that if $\chi$ is a linear
combination 
of $\chi_N$'s, then $$ 4 \cdot \chi((2,2)\mbox{-cycle}) = 3 \cdot \chi(
3\mbox{-cycle}) + \chi((3,2,2)\mbox{-cycle}). $$  Alternately, it
tells us 
that for any $N \ge 0$, four times the number of nonnegative integral 
solutions 
to $2 x_1 + 2 x_2 + x_3 + x_4 + x_5 = N$ is equal to three times the 
number of
such solutions to $3 x_1 + x_2 + x_3 + x_4 + x_5 = N$ plus the number
of such 
solutions to $3 x_1 + 2 x_2 + 2 x_3 = N$.
\end{example}

\section{Lower Bounds on $D(n)$}

Let $\lambda=(\lambda_1,\dots,\lambda_k)\vdash n$. The rational
function $f_\lambda(q)$ of equation (\ref{eq:f}) can be written as
  $$ f_\lambda(q) = p_\lambda(1,q,q^2,\dots), $$
where $p_\lambda$ denotes a power sum symmetric function. (See
\cite[Ch.\ 7]{EC2} for the necessary background on symmetric
functions.) Since the $p_\lambda$ for $\lambda\vdash n$ form a basis
for the vector space (say over $\C$) $\Lambda^n$ of all
homogeneous symmetric functions of degree $n$ \cite[Cor.\ 7.7.2]{EC2},
it follows that if $\{u_\lambda\}_{\lambda\vdash n}$ is any basis for
$\Lambda^n$ then 
  $$ D(n) = \dim\mathrm{span}_{\C}\{
  u_\lambda(1,q,q^2,\dots)\,:\,\lambda\vdash n\}. $$
In particular, let $u_\lambda=e_\lambda$, the elementary symmetric
function indexed by $\lambda$. Define
  $$ d(\lambda) = \sum_i {\lambda_i\choose 2}. $$
According to \cite[Prop.\ 7.8.3]{EC2}, we have
  $$ e_\lambda(1,q,q^2,\dots)=\frac{q^{d(\lambda)}}
    {\prod_i (1-q)(1-q^2)\cdots (1-q^{\lambda_i})}. $$
Since power series of different degrees (where the \emph{degree} of a
power series is the exponent of its first nonzero term) are
linearly independent, we obtain from Proposition~\ref{prop:dfn} the
following result. 

\begin{prop}
  Let $E(n)$ denote the number of distinct integers $d(\lambda)$,
  where $\lambda$ ranges over all partitions of $n$. Then $D(n)\geq
  E(n)$.
\end{prop}

\textsc{Note.} We could also use the basis $s_\lambda$ of Schur
functions instead of $e_\lambda$, since by \cite[Cor.\ 7.21.3]{EC2}
the degree of the power series $s_\lambda(1,q,q^2,\dots)$ is
$d(\lambda')$, where $\lambda'$ denotes the conjugate partition to
$\lambda$. 

Define $G(n)+1$ to be the least positive integer that cannot be
written in the form $\sum_i {\lambda_i\choose 2}$, where $\lambda
\vdash n$. Thus all integers $1,2,\dots,G(n)$ can be so represented,
so $D(n)\geq E(n)\geq G(n)$.  We can obtain a relatively tractable
lower bound for $G(n)$, as follows.  For a positive integer $m$, write
(uniquely)
  \begin{equation} m={k_1\choose 2}+{k_2\choose 2}+\cdots+{k_r\choose 2}, 
   \label{eq:nk2} \end{equation}
where $k_1\geq k_2\geq\cdots\geq k_r\geq 2$ and $k_1, k_2,\dots$ are
chosen successively as large as possible so that
  $$ m-{k_1\choose 2}-{k_2\choose 2}-\cdots -{k_i\choose 2}\geq 0 $$
for all $1\leq i\leq r$. For instance, $26={7\choose 2}+{3\choose 2}
+{2\choose 2}+{2\choose 2}$. Define $\nu(m)=k_1+k_2+\cdots+k_r$. Suppose
that $\nu(m)\leq n$ for all $m\leq N$. Then if $m\leq N$ we can write $m
={k_1\choose 2}+\cdots+{k_r\choose 2}$ so that $k_1+\cdots+k_r\leq
n$. Hence if $\lambda=\left(k_1,\dots,k_r,1^{n-\sum k_i}\right)$ (where
$1^s$ denotes $s$ parts equal to 1), then $\lambda$ is a partition of
$n$ for which $\sum_i {\lambda_i\choose 2}= m$. It follows that 
if $\nu(m)\leq n$ for all $m\leq N$ then $G(n)\geq N$. Hence if we
define $H(n)$ to be the largest integer $N$ for which  $\nu(m)\leq n$
whenever $m \leq N$, then we have established the string of
inequalities
  \begin{equation} D(n) \geq E(n) \geq G(n) \geq H(n). 
    \label{eq:string} \end{equation}
Here is a table of values of these numbers for $1\leq n\leq 23$. Note
that $D(n)$ appears to be close to $E(n+1)$. We don't have any
theoretical explanation of this observation.

\begin{table}[h]
\begin{tabular}{|c|r|r|r|r|r|r|r|r|r|r|r|r|r|r|}
\hline
 $n$   & 1 & 2 & 3 & 4 & 5 & 6 & 7 & 8 & 9 & 10 & 11 & 12 & 13 & 14 \\ 
\hline
$D(n)$ & 1 & 2 & 3 & 5 & 7 & 11 & 13 & 19 & 23 & 29 & 35 & 45 & 51 & 62  \\
$E(n)$ & 1 & 2 & 3 & 5 & 7 &  9 & 13 & 18 & 21 & 27 & 34 & 39 & 46 & 54  \\ 
$G(n)$ & 0 & 1 & 1 & 3 & 4 & 4 & 7 & 13 & 13 & 18 & 25 & 32 & 32 & 32  \\
$H(n)$ & 0 & 1 & 1 & 3 & 4 & 4 & 7 & 11 & 13 & 18 & 19 & 19 & 25 & 32  
 \\ \hline \end{tabular}

\vskip 0.1in

\begin{tabular}{|c|r|r|r|r|r|r|r|r|r|}
\hline
 $n$ &  15 & 16 & 17 & 18 & 19 & 20 & 21 & 22 & 23  \\ \hline
$D(n)$  & 69 & 79 & 90 & 106 & 118 & 134 & 146 & 161 & 176  \\
$E(n)$  & 61 & 72 & 83 & 92 & 106 & 118 & 130 & 145 & 162  \\ 
$G(n)$  & 40 & 49 & 52 & 62 & 73 & 85 & 102 & 112 & 127\\
$H(n)$ & 40 & 43 & 52 & 62 & 73 & 85 & 89 & 102 & 116  \\ \hline
\end{tabular}

\vskip 0.1in

\caption{Values of $D(n)$, $E(n)$, $G(n)$, $H(n)$ for small $n$}
\end{table}

\begin{prop} \label{prop:num}
We have
  \begin{equation} \nu(m) \leq \sqrt{2m}+3m^{1/4} 
      \label{eq:nun} \end{equation}
for all $m\geq 405$.
\end{prop}

\begin{proof}
The proof is by induction on $m$. It can be checked with a computer
that equation (\ref{eq:nun}) is true for $405\leq m\leq 50000$. Now
assume that $M>50000$ and that (\ref{eq:nun}) holds for $405\leq m
<M$. Let $p=p_M$ be the unique positive integer
satisfying
  $$ {p\choose 2}\leq M < {p+1\choose 2}. $$
Thus $p$ is just the integer $k_1$ of equation (\ref{eq:nk2}).
Explicitly we have 
  $$ p_M = \left\lfloor\frac{1+\sqrt{8M+1}}{2}\right\rfloor. $$
By the definition of $\nu(M)$ we have
  $$ \nu(M) = p_M + \nu\left( M-{p_M\choose 2}\right). $$
It can be checked that the maximum value of $\nu(m)$ for $m<405$ is
$\nu(404) =42$. Set $q_M=(1+\sqrt{8M+1})/2$. Since $M-{p_M\choose 2}
\leq p_M\leq q_M$, by the induction hypothesis we have 
  $$ \nu(M) \leq q_M +\max(42, \sqrt{2q_M}+ 3q_M^{1/4}). $$
It is routine to check that when $M>50000$ the right hand side is less
than $\sqrt{2M}+3M^{1/4}$, and the proof follows.
\end{proof}

\begin{prop} \label{prop:c}
There exists a constant $c>0$ such that
  $$ H(n) \geq \frac{n^2}{2}-cn^{3/2} $$
for all $n\geq 1$.
\end{prop}

\begin{proof}
{}From the definition of $H(n)$ and Proposition~\ref{prop:num} (and the
fact
that the right-hand side of equation (\ref{eq:nun}) is increasing), 
along with the
inquality $\nu(m) \le 42 = \lceil\sqrt{2 \cdot 405}+3 \cdot
405^{1/4}\rceil$ for $m \le 404$, it follows that 
  $$ H\left(\lceil\sqrt{2m}+3m^{1/4}\rceil\right) \geq m $$
for $m>404$.  For $n$ sufficiently large, we can evidently choose $m$
such that $n = \lceil\sqrt{2m}+3m^{1/4}\rceil$, so 
$H(n) \geq m$.  Since $ \sqrt{2m} + 3m^{1/4} + 1 > n$, an application
of the quadratic formula (again for $n$ sufficiently large) shows
$$ m^{1/4} \geq \frac{-3 + \sqrt{9 + 4\sqrt{2} (n-1)}}{2 \sqrt{2}}, $$
from which the result follows without difficulty.

\end{proof}

Since we have established both upper bounds (equation (\ref{eq:dn}))
and lower bounds (equation (\ref{eq:string}) and
Proposition~\ref{prop:c}) for $D(n)$ asymptotic to $n^2/2$, we obtain
the following corollary.

\begin{cor}
There holds the asymptotic formula $D(n)\sim \frac 12 n^2$.
\end{cor}

\end{document}